\documentclass[11pt]{article}
\usepackage[utf8]{inputenc}
\usepackage[english]{babel}
\usepackage{amsmath}
\usepackage{amsfonts}
\usepackage{amssymb}
\usepackage{epsfig}
\usepackage{graphics,psfrag,rotating}
\usepackage{graphicx}
\usepackage{dcolumn}
\usepackage{bm}
\usepackage{color}
\usepackage[usenames,dvipsnames,svgnames]{xcolor}
\usepackage[T1]{fontenc}
\usepackage{multirow}
\usepackage{float}

\usepackage{subfigure}
\usepackage{pb-diagram} %=============for diagrams

\usepackage{enumitem}
\usepackage[colorlinks=true,
            linkcolor=red,
          urlcolor=gray,
            citecolor=blue]{hyperref}

\newtheorem{Def}{Definition}
\newtheorem{Th}{Theorem}
\newtheorem{Prop}{Proposition}
\newtheorem{Lem}{Lemma}

\date{}
\begin{document}

% Название доклада
\title{\textbf{On the quotient of projective frame space and the Desargues theorem}}
\maketitle

% Информация об авторе
\begin{center}
\author{Artur Kuleshov$^{1}$\\
$^{1}$I. Kant Baltic Federal University, Nevskogo str. 14, 236016
Kaliningrad, Russia}
\end{center}

\begin{abstract}
  We consider an $n$-dimensional projective space $\mathbb{P}_n$ ($n\geq2$) and a fixed point $A$ on it. Let $F(\mathbb{P}_n)$ be the manifold of all the projective frames of $\mathbb{P}_n$ having $A$ as their first vertice. We define the action of $G=St_{A}\subset GP(n)$ on $F(\mathbb{P}_n)$ in a natural way. The Lie group epimorphism $\beta\colon G\to GL(V)$ acts as follows $g\mapsto d_{A}g$ where $V=T_{A}\mathbb{P}_n$. We study the geometry of orbit space $\Phi(\mathbb{P}_n)$ of projective frame space $F(\mathbb{P}_n)$ under the action of the kernel $H$ of this epimorphism $\beta$. By applying some $n$-dimensional version of the Desargues theorem we could get a purely geometrical description of such $H$-orbits.

%We study the gravitational collapse of homogeneous and   %%%%%%%%%%%%%%% We study (...) in the framework of (...)
%in-homoge\-neous perfect fluid in the framework of general %%%%%%%%%%%%% By (...) we could get (...)
%relativity. By adding a suitable Lagrangian in the
%Einstein-Hilbert action to the already existing Ricci scalar 'R',
%we could get the freedom to explain cosmic acceleration. In this
%new domain of general relativity in presence of dark energy, which
%in this case is the one corresponding to Lagrangian $-X-V$ (where
%$X$ and $V$ are the kinetic and potential energy of the field
%respectively), called the phantom fluid, the gravitational
%collapse of ordinary matter perfect fluid is investigated.
%Formulation of the scalar field corresponding to the linear
%phantom fluid having the above mentioned Lagrangian has been done.
\end{abstract}

%======================

\section{Introduction}

According to the prolongations and scopes method \cite{Laptev_1953} of studying geometrical structures on submanifolds immersed into homogeneous spaces ($\mathbb{P}_n$ in our case) frame bundles over such submanifolds and quotients of these bundles are considered (see, e.g., \cite{Akivis_Goldberg_1993}, \cite{Belova_2008}, \cite{Polyakova_2008}, \cite{Schevchenko}). Such principal bundles and their structure groups are studied purely in terms of structure equations of their Maurer~-- Cartan forms \cite{Ivey_landsberg} because it is sufficient for obtaining local results considered in these papers. But there is a natural question about the explicit description of such quotient bundles in purely geometrical terms.

In the paper we restrict our attention to one of the simplest cases. Namely, let's consider a projective $n$-dimensional space $\mathbb{P}_n$ together with its distinguished point $A\in \mathbb{P}_n$. Let $G$ be the stabilizer of $A$ in the projective transformation group of $\mathbb{P}_n$. The structure equations of $G$ are the following
$$
  d\omega^i_j=\omega^k_j\wedge\omega^i_k,\quad d\omega_i=\omega^k_i\wedge\omega_k,
$$
where $\omega^i_j$, $\omega_i$ are the Maurer~-- Cartan forms of $G$. From these equations it follows that there exists a quotient group $\bar{G}$ of $G$ isomorphic (at least locally) to $GL(n)$. But it is not easy to describe its action on geometric images (i.e. points, frames etc.). In particular one can ask the following question: \textit{let $H$ be the corresponding normal subgroup, how to describe $H$-orbits of projective frames?} Its importance comes from the following facts. Firstly, as we will see further, these $H$-orbits can be identified with bases of some vector space. And, secondly, any action of elements of $G$ sends $H$-orbit as a whole object to another $H$-orbit.

\smallskip

%=====================================================================================

\section{Basic concepts and claims}

According to the approach proposed in \cite{Berger} we start with the following definition.

\begin{Def}\label{projecive_space}
  Let $V_{n+1}$ be $(n+1)$-dimensional vector space and let $\sim$ be the collinearity relation on $V^0_{n+1}=V_{n+1}\backslash\{\vec{0}\}$. Then an $n$-dimensional projective space is  $\mathbb{P}_n$ is the quotient space $\mathbb{P}_n=V^0_{n+1}/{\sim}$. The vector space $V_{n+1}$ is said to be associated with $\mathbb{P}_n$. The canonical surjection $\pi\colon V^0_{n+1}\to \mathbb{P}_n$ sends each vector $\vec{a}\in V^0_{n+1}$ to its equivalence class $[\vec{a}]$.
\end{Def}

\begin{Def}\label{projecive_frame}
  A projective frame in $\mathbb{P}_n$ is an ordered set $\mathcal{R}$ of $n+2$ points $A_0,\,A_1,\ldots,\,A_n,\,E\in\mathbb{P}_n$ such that any $n+1$ points of $\mathcal{R}$ are in generic position:
  \begin{equation}\label{proj_frame}
    \mathcal{R}=\{A_0,\,A_1,\ldots,\,A_n,\,E\}.
  \end{equation}
\end{Def}

\begin{Def}
  Say a basis
  $$
    \vec{\mathcal{R}}=\{\vec{A_0},\,\vec{A_1},\ldots,\,\vec{A_n}\}
  $$
  of $V_{n+1}$ generates the frame (\ref{proj_frame}) iff
  \begin{gather*}
    \pi(\vec{A}_0)=A_0,\,\pi(\vec{A}_1)=A_1,\,\ldots,\,\pi(\vec{A}_n)=A_n,\\
    \pi(\vec{A}_0 + \vec{A}_1+ \ldots + \vec{A}_n)=E.
  \end{gather*}
\end{Def}

Let us distinguish some point $A\in \mathbb{P}_n$ and call it a {\it center} of $\mathbb{P}_n$.

\begin{Def}
  We say a projective frame $\mathcal{R}=\{A_0,\,A_1,\ldots,\,A_n,\,E\}$ to be adapted (or, equivalently, centroprojective) if its first vertice $A_0$ coinsides with the center, i.e. $A_0=A$.
\end{Def}

\begin{Def}
  Homogeneous coordinates of a point $M$ w.r.t. a frame $\mathcal{R}$ are the coordinates of some vector $\vec{M}\in\pi^{-1}(M)$ w.r.t. a basis $\vec{\mathcal{R}}$ generating the frame $\mathcal{R}$:
  $$
    M(x^0\colon\,x^1\colon\,\ldots\,\colon x^n)_{\mathcal{R}}\quad \Leftrightarrow \quad\vec{M}=x^0\vec{A}_0 + x^1\vec{A}_1 + \ldots\ + x^n\vec{A}_n.
  $$
  By $U_{\mathcal{R}}$ denote the open subset of $\mathbb{P}_n$ given by $x^0\neq0$. Affine (non-homogeneous) coordinates of the point $M\in U_{\mathcal{R}}$ w.r.t. the frame $\mathcal{R}$ are the following quotients:
  $$
    X^i=\frac{x^i}{x^0},\quad i=\overline{1,n}.
  $$
\end{Def}

{\bf Remark}. The point $A$ is given by $(x^0:\,0:\,\ldots:\,0)$, $x^0\neq0$, for any adapted frame ${\mathcal{R}}$, so $U_{\mathcal{R}}$ is an open neighborhood of $A$.

\smallskip

\begin{Def}
  By $\varphi_{\mathcal{R}}$ we denote the affine chart on $U_{\mathcal{R}}$, i.e. the mapping sending each point $M\in U_{\mathcal{R}}$ to its non-homogeneous coordinates:
  $$
    \varphi_{\mathcal{R}}\colon U_{\mathcal{R}}\to \mathbb{R}^n,\quad M\mapsto (X^1\,\ldots,\,X^n).
  $$
\end{Def}

\begin{Prop}\label{Change_frames_formulas}
  Let $\mathcal{R}$ and $\mathcal{R}'$ be some frames such that $\mathcal{R}$ is adapted. Then the following conditions are equivalent:

  1)~$\mathcal{R}'$ is adapted also;

  2)~for any bases ${\vec{\mathcal{R}}}$ and ${\vec{\mathcal{R}}'}$ generating these frames there exist coefficients $a^0_0,\,a^0_i,\,a^i_j$ ($i,\,j=\overline{1,n}$) such that
  \begin{equation*}
    \vec{A}'_0 = a^0_0 \vec{A}_0,\quad \vec{A}'_i = a^0_i \vec{A}_0 + a^j_i \vec{A}_j,\quad a^0_0\neq0,\quad\det(a^j_i)\neq0;
  \end{equation*}

  3)~there exist unique coefficients $\alpha_i$, $\alpha^i_j$ such that for any point $M$ its affine coordinates change under the law
  \begin{equation}\label{transition_equations}
    {X}^i=\frac{\alpha^i_j\tilde{X}^j}{1+\alpha_j \tilde{X}^j}.
  \end{equation}
\end{Prop}

\begin{Def}
  The equations (\ref{transition_equations}) are called transition equations of a pair $(\mathcal{R},\,\mathcal{R}')$.
\end{Def}

\begin{Def}\label{projective_transformation}
  The mapping $f\colon\mathbb{P}_n\to \mathbb{P}_n$ is called a projective transformation $\Leftrightarrow$ there exists a non-degenerate linear operator $u\colon V_{n+1}\to V_{n+1}$ such that the following diagram commutes:
  $$
  \begin{diagram}
    \node{V^0_{n+1}}
      \arrow{e,t}{u}   %===============arrow to the left (----->)
      \arrow{s,l}{\pi}   %===============arrow to the down (|)
    \node{V^0_{n+1}}
      \arrow{s,l}{\pi}\\  %===============arrow to the down (|)
    \node{\mathbb{P}_n}
      \arrow{e,r}{f} %===============arrow to the left (----->)
    \node{\mathbb{P}_n}
  \end{diagram}
  $$
\end{Def}

{\bf Remark}. $u$ is determined by $f$ up to a non-zero scalar multiple.

\smallskip

We introduce the following notation:

\smallskip

\begin{center}
\begin{tabular}{|rl|}
  \hline
  % после \\: \hline или \cline{col1-col2} \cline{col3-col4} ...
  $F(P_n)$ & --- the set of all the adapted projective frames \\
  $V$ & --- the tangent space $T_A(\mathbb{P}_n)$ to $\mathbb{P}_n$ at $A$ \\
  $F(V)$ & --- the set of all the linear frames (i.~e. bases) of $V$ \\
  $GP(n)$ & --- the projective transformation group  of $\mathbb{P}_n$ \\
  $G$ & --- the stabilizer of $A$ in $GP(n)$ \\
  $GL(V)$ & --- the group of non-degenerate linear operators in $V$ \\
  \hline
\end{tabular}
\end{center}

\smallskip

Let $\mathcal{R}\in F(P_n)$ and $\varepsilon\in F(V)$, where $$
\mathcal{R}=\{A_0,\,A_1,\ldots,\,A_n,\,E\},\quad \varepsilon=\{\vec{e}_1,\ldots,\vec{e}_n\}.
$$
We denote
\begin{gather*}
  g\cdot\mathcal{R}\colon=\{g(A_0),\,g(A_1),\ldots,g(A_n),\,g(E)\},\quad g\in G;\\ \psi\cdot\varepsilon\colon=\{\psi(\vec{e}_1),\ldots,\psi(\vec{e}_n)\},\quad \psi\in GL(V).
\end{gather*}

\begin{Prop}\label{frames_and_actions}
  1)~$g\cdot\mathcal{R}\in F(\mathbb{P}_n)$ for any $g\in G$, $\mathcal{R}\in F(P_n)$.

  2)~For any two frames $\mathcal{R},\,\mathcal{R}'\in F(\mathbb{P}_n)$ there exists a unique $g\in G$ such that $\mathcal{R}'=g\cdot \mathcal{R}$.
\end{Prop}

\begin{Def}
  Let $g\in G$, $\mathcal{R}\in F(\mathbb{P}_n)$. Equations of $g$ w.r.t. $\mathcal{R}$ are the transtion equations of the pair $(\mathcal{R},\,g\cdot \mathcal{R})$.
\end{Def}

\smallskip

\begin{Def}
  A frame action of ${G}$ on $F(P_n)$ is the action $$q_P\colon {G}\times F(P_n)\to F(P_n)$$ defined as follows:
  $$
    (f, \mathcal{R})\mapsto f\cdot\mathcal{R},\quad f\in G,\,\, \mathcal{R}\in F(P_n).
  $$
\end{Def}

\smallskip

\begin{Def}
  A frame action of $GL(V)$ on $F(V)$ is the action $$q_V\colon {G}\times F(V)\to F(V)$$ defined as follows:
  $$(\psi, \varepsilon)\mapsto \psi\cdot\varepsilon,\quad f\in G,\,\, \mathcal{R}\in F(P_n).$$
\end{Def}

\smallskip

\begin{Prop}
  The mappings $q_P$ and $q_V$ are well-defined free transitive smooth actions of the Lie groups $G$ and $GL(V)$ on $F(P_n)$ and $F(V)$ respectively.
\end{Prop}

\smallskip

%=====================================================================================

\section{$H$-orbits and quotient frame space}

\begin{Def}
  A linearizing mapping is the mapping $\beta\colon G\to GL(V)$ acting as follows
  $$\beta\colon f\mapsto d_A f.$$
\end{Def}

\begin{Lem} (see \cite{Conlon})
  Let $\varphi\colon G\to G'$ be a continuous homomorphism of Lie groups and let $H=\ker\varphi$. Then $H$ is a properly embedded, normal Lie subgroup of $G$, $G/H$ is canonically a Lie group, and the induced map $\bar{\varphi}\colon G/H\to G'$ is an injective immersion of this Lie group as a Lie subgroup of $G'$.
\end{Lem}

\begin{Prop}

  1)~$\beta$ is a Lie group epimorphism.

  2)~The kernel $H:=\ker \beta$ is a closed normal subgroup of $G$.

  3)~$\bar{\beta}\colon fH\mapsto d_A f$ is a canonical Lie group isomorphism between $G/H$ and $GL(V)$.
\end{Prop}

\smallskip

%{\it Proof}. 1)~follows from the chain rule for the differential: $$d_A(f\circ g)=d_A{f}\circ d_A{g},\quad f,\,g\in G.$$ 2) and 3) follow from 1) and Lemma~1. $\square$
{\bf Remark}. This proposition allows us not to distinguish the Lie groups $G/H$ and $GL(V)$.

\smallskip

\begin{Def}
  We say the group $\bar{G}:=G/H=GL(V)$ is a linear quotient group of $G$.
\end{Def}

\smallskip

\begin{Prop}\label{coefficients__of_h}
  Let $g\in G$, $\mathcal{R}\in F(\mathbb{P}_n)$, and let $(\ref{transition_equations})$ be the equations of $g$ w.r.t. $\mathcal{R}$. Then
  $$
    g\in H\quad\Leftrightarrow\quad \alpha^i_j=\delta^i_j\quad \Leftrightarrow \quad {X}^i=\frac{\tilde{X}^i}{1+\alpha_j \tilde{X}^j}.
  $$
\end{Prop}

\begin{Def}\label{equiv_frames}
  Two frames $\mathcal{R}$ and $\mathcal{R}'$ are said to be equivalent ($\mathcal{R}\sim\mathcal{R}'$) if they belong to the same $H$-orbit. We denote by $[\mathcal{R}]$ the equivalence class of the frame $\mathcal{R}$, by $\Phi(P_n)$ the set of all such equivalence classes, and by $p\colon F(\mathbb{P}_n)\to \Phi(\mathbb{P}_n)$ denote the canonical projection $\mathcal{R}\mapsto [\mathcal{R}]$.
\end{Def}

\begin{Prop}\label{p}
  $p$ is a surjective submersion.
\end{Prop}

\begin{Prop}
  Every action of $G$ on $F(\mathbb{P}_n)$ sends $H$-orbits to $H$-orbits, i.e.
  $$
    (\forall f\in G)(\forall\,\mathcal{R},\,\mathcal{R}'\in F(\mathbb{P}_n))(\mathcal{R}\sim\mathcal{R}'\,\Rightarrow f\cdot \mathcal{R}\sim f\cdot \mathcal{R}').
  $$
\end{Prop}

\smallskip

\begin{Def}
  A quotient action of $G$ on $\Phi(P_n)$ is the action $$\gamma\colon {G}\times\Phi(P_n)\to \Phi(P_n)$$ defined as follows:
  $$f\cdot [\mathcal{R}]\colon=[f\cdot \mathcal{R}],\quad f\in G,\,\, \mathcal{R}\in F(P_n).$$
\end{Def}

\smallskip

\begin{Prop}
  Any two elements $f,\,g\in G$ belonging to the same $H$-coset act on $\Phi(\mathbb{P}_n)$ by the same way i.e.
  $$
    (\forall f,\,g\in G)(\forall \mathcal{R}\in F(\mathbb{P}_n))(f\in gH\,\Rightarrow f\cdot [\mathcal{R}]=g\cdot [\mathcal{R}]).
  $$
\end{Prop}

\smallskip

\begin{Def}\label{bar_gamma}
  A quotient action of $\bar{G}$ on $\Phi(P_n)$ is the action $\bar{\gamma}\colon \bar{G}\times\Phi(P_n)\to \Phi(P_n)$ defined as follows:
  $$gH\cdot [\mathcal{R}]\colon=[g\cdot \mathcal{R}],\quad g\in G,\,\, \mathcal{R}\in F(P_n).$$
\end{Def}

\begin{Prop}
  The following diagram commutes for any $g\in G$:
  $$
  \begin{diagram}
    \node{F(\mathbb{P}_n)}
      \arrow{e,t}{g}   %===============arrow to the left (----->)
      \arrow{s,l}{p}
    \node{F(\mathbb{P}_n)}
      \arrow{s,l}{p}\\
    \node{\Phi(\mathbb{P}_n)}
      \arrow{e,r}{gH}
    \node{\Phi(\mathbb{P}_n)}
  \end{diagram}
  $$
\end{Prop}

{\it Proof} follows immediately from Definition~\ref{bar_gamma}. $\square$

\begin{Prop}
  $\bar{\gamma}$ is free, transitive and smooth.
\end{Prop}

\smallskip
%=====================================================================================

\section{Isomorphism of $\bar{G}$-spaces $\Phi(\mathbb{P}_n)$ and $F(V)$}

\begin{Def}\label{alpha}
  Let $\alpha\colon F(\mathbb{P}_n)\to F(V)$ be the map acting as follows:
  $$
    \alpha\colon \mathcal{R}\mapsto\varepsilon(\varphi_{\mathcal{R}})
  $$
  where $\varepsilon(\varphi_{\mathcal{R}})$ is the natural basis of the space $V$ generated by $\varphi_{\mathcal{R}}$.
\end{Def}

{\bf Remark}. The mapping $\alpha$ is a surjective submersion.

\smallskip

\begin{Prop}
  The following diagram commutes for any $g\in G$:
  $$
  \begin{diagram}
    \node{F(\mathbb{P}_n)}
      \arrow{e,t}{g}   %===============arrow to the left (----->)
      \arrow{s,l}{\alpha}
    \node{F(\mathbb{P}_n)}
      \arrow{s,l}{\alpha}\\
    \node{F(V)}
      \arrow{e,r}{d_A{g}}
    \node{F(V)}
  \end{diagram}
  $$
\end{Prop}

\begin{Prop}
  Equivalence classes of frames are exactly preimages under the map $\alpha$, i.e.
  $$
    (\forall \mathcal{R},\,\mathcal{R}'\in F(\mathbb{P}_n))(\mathcal{R}\sim\mathcal{R}'\,\Leftrightarrow \alpha(\mathcal{R})=\alpha(\mathcal{R}')).
  $$
\end{Prop}

\begin{Def}\label{bar_alpha}
  Let $\bar{\alpha}\colon \Phi(\mathbb{P}_n)\to F(V)$ be the map acting as follows: $$\bar{\alpha}\colon [\mathcal{R}]\mapsto\varepsilon(\varphi_{\mathcal{R}}).$$
\end{Def}

\begin{Prop}
  The following diagram commutes:
  $$
  \begin{diagram}
    \node{F(\mathbb{P}_n)}
      \arrow{e,t}{\alpha}   %===============arrow to the left (----->)
      \arrow{s,l}{p}
    \node{F(V)}\\
    \node{\Phi(\mathbb{P}_n)}
      \arrow{ne,r}{\bar{\alpha}}
  \end{diagram}
  $$
\end{Prop}

{\it Proof} follows immediately from Definitions~\ref{equiv_frames}, \ref{alpha} and \ref{bar_alpha}. $\square$

\begin{Prop}\label{feature_of_bar_alpha}
  $\bar{\alpha}$ is a diffeomorphism.
\end{Prop}

\begin{Prop}\label{main_diagram}
  The following diagram commutes for any $g\in G$:
  $$
  \begin{diagram}
    \node{\Phi(\mathbb{P}_n)}
      \arrow[3]{e,t}{gH}   %===============arrow to the left (----->)
      \arrow[2]{s,l}{\bar{\alpha}}
    \node[3]{\Phi(\mathbb{P}_n)}
      \arrow[2]{s,r}{\bar{\alpha}}\\
    \node[2]{F(\mathbb{P}_n)}
      \arrow{nw,t}{p}
      \arrow{e,t}{g}    %===============arrow to the left (----->)
      \arrow{sw,b}{\alpha}
    \node{F(\mathbb{P}_n)}
      \arrow{ne,t}{p}
      \arrow{se,b}{\alpha}\\
    \node{F(V)}
      \arrow[3]{e,b}{d_A{g}}   %===============arrow to the left (----->)
    \node[3]{F(V)}
  \end{diagram}
  $$
\end{Prop}

\begin{Th}
  $\bar{\alpha}\colon \Phi(\mathbb{P}_n)\to F(V)$ is an isomorphism of $\bar{G}$-spaces.
\end{Th}
{\it Proof} follows immediately from Propositions~\ref{main_diagram}, \ref{p} and \ref{feature_of_bar_alpha}. $\square$

\smallskip
%=====================================================================================

\section{Perspectivity and the Desargues theorem}

Further on we restrict ourselves to the case $n\geq2$. Let $\mathcal{R},\,\mathcal{R}'\in F(P_n)$, where
\begin{equation*}
  \mathcal{R}=\{A_0,\,A_1,\ldots,\,A_n,\,E\},\quad \mathcal{R}'=\{A'_0,\,A'_1,\ldots,\,A'_n,\,E'\}.
\end{equation*}
Recall that $A_0 = A'_0 = A$. Consider two bases $\vec{\mathcal{R}}$ and $\vec{\mathcal{R}}'$ generating the frames $\mathcal{R}$ and $\mathcal{R}'$ respectively:
\begin{equation}\label{eq1}
  \vec{\mathcal{R}}=\{\vec{A_0},\,\vec{A_1},\ldots,\,\vec{A_n}\},\quad    \vec{\mathcal{R}}'=\{\vec{A}'_0,\,\vec{A}'_1,\ldots,\,\vec{A}'_n\}.
\end{equation}

\smallskip

\begin{Def}
  Two frames $\mathcal{R}$ and $\mathcal{R}'$ are said to be in perspective if
  $$A'_i\in A_i A_0\quad (i=\overline{1,n}),\quad E'\in E A_0.$$
\end{Def}

\smallskip

\begin{Prop}\label{Perspectivity}
  Let $\mathcal{R}$ and $\mathcal{R}'$ be any adapted frames. Then the following conditions are equivalent:

  1)~$\mathcal{R}$ and $\mathcal{R}'$ are in perspective;

  2)~for any bases ${\vec{\mathcal{R}}}$ and ${\vec{\mathcal{R}}'}$ generating these frames there exist coefficients $b^0_0,\,b^0_i,\,c^0_0$ such that
  \begin{equation}\label{eq*a}
    \vec{A}'_0 = b^0_0 \vec{A}_0,\quad \vec{A}'_i = b^0_i \vec{A}_0 + c^0_0 \vec{A}_i,\quad b^0_0\neq0,\,c^0_0\neq0,\quad i=\overline{1,n}.
  \end{equation}

  3)~for any basis $\vec{\mathcal{R}}$ generating $\mathcal{R}$ there are unique numbers $a_1,\,\ldots,\,a_n$, $h\neq0$ and a unique basis $\vec{\mathcal{R}}''$ generating $\mathcal{R}'$ such that the following equalities are hold:
  \begin{equation}\label{eq*4}
    \vec{A}''_0 = \vec{A}_0,\quad \vec{A}''_i = a_i\vec{A}_0+ h\vec{A}_i,\quad i=\overline{1,n}.
  \end{equation}

  4)~there exist unique coefficients $h\neq0$, $a_1,\,\ldots,\,a_n$ such that for any point $M$ its affine coordinates change under the law
  \begin{equation}\label{eq_p}
    {X}^i=\frac{h \tilde{X}^i}{1+a_j \tilde{X}^j}.
  \end{equation}
\end{Prop}

%{\it Proof}. 2)~$\Rightarrow$~3). We prove the existence. Let $\vec{\mathcal{R}}$ and $\vec{\mathcal{R}}''$ be any bases generating the frames. Due to 2) there exist coefficients $b^0_0,\,b^0_i,\,c^0_0$ such that (\ref{eq*a}) hold. Let
%$$
%  \vec{A}''_0 = \frac{1}{b^0_0}\vec{A}'_0,\quad \vec{A}''_i = \frac{1}{b^0_0}\vec{A}'_i,\quad a_i=\frac{b^0_i}{b^0_0},\quad h=\frac{c^0_0}{b^0_0}.
%$$
%Then (\ref{eq*4}) are hold. The uniqueness part is obvious.
%
%2)~$\Rightarrow$~3) is obvious.
%
%($\square$)

%{\bf Remark~1}. The statement 2) is equivalent to its weak version:
%  ~{\it there \underline{exist} bases (\ref{eq1}) generating the frames $\mathcal{R}$ and $\mathcal{R}'$ such that for some numbers $a_1,\,\ldots,\,a_n$, $h$ (\ref{eq*4}) hold.}

%{\bf Remark~2}. The equalities (\ref{eq*4}) imply that $\vec{E}'' = e\vec{A}_0+ \vec{E}$, where $e:=a_1+\ldots+a_n+1-h$.

\smallskip

\begin{Def}
  Transform coefficients of a pair $(\mathcal{R},\,\mathcal{R}')$ of frames in perspective are the numbers $a_1,\,\ldots,\,a_n$, $h$ determined in Proposition~\ref{Perspectivity}.
\end{Def}

{\bf Remark~1}. Proposition~\ref{Perspectivity} implies that the transform coefficients don't depend on the choice of $\vec{\mathcal{R}}$. So they are completely determined by the pair $(\mathcal{R},\,\mathcal{R}')$.

{\bf Remark~2}. In a particular case $a_i=0$ for all $i=\overline{1,n}$ the transformation (\ref{eq_p}) is homothetic in the affine chart $(\varphi_{\mathcal{R}},\,U_{\mathcal{R}})$.

\smallskip

\begin{Prop}\label{equivalence}
  Let $\mathcal{R}$ and $\mathcal{R}'$ be any adapted frames. Then the following conditions are equivalent:

  1)~$\mathcal{R}\sim\mathcal{R}'$;

  2)~for any bases ${\vec{\mathcal{R}}}$ and ${\vec{\mathcal{R}}'}$ generating these frames there exist coefficients $a^0_0,\,a^0_i$ such that
  \begin{equation*}
    \vec{A}_0 = a^0_0 \vec{A}'_0,\quad \vec{A}_i = a^0_i \vec{A}'_0 + a^0_0 \vec{A}'_i,\quad a^0_0\neq0.
  \end{equation*}

  3)~for any basis $\vec{\mathcal{R}}$ generating $\mathcal{R}$ there are unique numbers $a_1,\,\ldots,\,a_n$, and a unique basis $\vec{\mathcal{R}}''$ generating $\mathcal{R}'$ such that the following equalities are hold:
  \begin{equation*}
    \vec{A}''_0 = \vec{A}_0,\quad \vec{A}''_i = a_i\vec{A}_0+ \vec{A}_i,\quad i=\overline{1,n};
  \end{equation*}

  4)~there exist unique coefficients $a_1,\,\ldots,\,a_n$ such that for any point $M$ its affine coordinates change under the law
  \begin{equation}\label{eq_q}
    {X}^i=\frac{\tilde{X}^i}{1+a_j \tilde{X}^j}.
  \end{equation}
\end{Prop}

\begin{Th}\label{Criterion_1}
  $\mathcal{R}\sim\mathcal{R}'$ $\Leftrightarrow$ $\mathcal{R}$ and $\mathcal{R}'$ are in perspective $\&$ the canonical coefficients of the pair $(\mathcal{R},\,\mathcal{R}')$ satisfy the following condition:
  \begin{equation}\label{eq3}
    h=1.
  \end{equation}
\end{Th}

{\it Proof} follows immediately from (\ref{eq_p}) and (\ref{eq_q}). $\square$

\smallskip

\begin{Def}
  We say the frames $\mathcal{R}$ and $\mathcal{R}'$ to be in strict perspective if they are in perspective and their corresponding points are not coincide, i.e.
  $$
    {A}'_i \neq {A}_i\,\, (i=\overline{1,n}),\quad E'\neq E.
  $$
\end{Def}

\smallskip

For any two frames $\mathcal{R}$ and $\mathcal{R}'$ in strict perspective the minimal subspace (with respect to inclusion) $\mathcal{L}_{(\mathcal{R},\,\mathcal{R}')}\subset P_n$ containing the set of points $B_{ij}$, $B_i$ is defined, where
$$
B_{ij}=A_i A_j\cap A'_i A'_j,\quad B_i=A_i E\cap A'_i E',\quad 1\leq i<j\leq n.
$$

\smallskip

\begin{Lem}\label{Desargues_classical} (the Desargues theorem, classical version)
  Let $A_1B_1$, $A_2B_2$ and $A_3B_3$ be three concurrent lines on $\mathbb{P}_2$. Then the points $$
    C_{12}=A_1A_2\cap B_1B_2,\quad C_{13}=A_1A_3\cap B_1B_3,\quad C_{23}=A_2A_3\cap B_2B_3
  $$
  are lying on a straight line.
\end{Lem}

{\it Proof} see, e.g., in \cite{Lelong-Ferrand}.$\square$

\begin{Lem}\label{Desargues_frames} (the Desargues theorem, frame version)
  Let $\mathcal{R},\,\mathcal{R}'\in \mathbb{P}_2$ be two frames in strict perspective. Then $\mathcal{L}_{(\mathcal{R},\,\mathcal{R}')}$ is a straight line.
\end{Lem}

%%%%%%%%%%%%%%%%%%%%%%%%%%%%%%%%%% Hyperplane theorem %%%%%%%%%%%%%%%%%%%%%%%%%%%%
\begin{Th}\label{Hyperplane}
  $\mathcal{L}_{(\mathcal{R},\,\mathcal{R}')}\subset P_n$ is a hyperplane in $\mathbb{P}_n$ for any two frames $\mathcal{R}$ and $\mathcal{R}'$ in strict perspective.
\end{Th}

\smallskip

{\it Proof} is based on \cite{Bell_1955}. In the case $n=2$ the statement is just Lemma~\ref{Desargues_frames}. Further on, we shall assume that $n>2$. For any collection of subsets $Y_1,\ldots,\, Y_s\subset \mathbb{P}_n$ we denote by $\langle Y_1,\ldots,\,Y_s\rangle$ the minimal plane in $\mathbb{P}_n$ containing all of them. Then for $1\leq i<j\leq n$ we have
\begin{gather*}
  B_{ij}\in A_i A_j \subset \langle A_1,\ldots A_n\rangle\,=:\mathcal{M},\\
  B_{ij}\in A'_i A'_j \subset \langle A'_1,\ldots A'_n\rangle\,=:\mathcal{M}'.
\end{gather*}
Let $\mathcal{N}=\mathcal{M}\cap\mathcal{M}'$. Then $B_{ij}\in \mathcal{N}$ and $\dim\mathcal{N}=n-2$ due to the conditions of the theorem. Consider 3-plane $L_{ij}:=\langle A_0 A_i,\,A_0 A_j,\,A_0 E\rangle$. The 2-planes $\langle A_i,\,A_j,\,E\rangle$ and $\langle A'_i,\,A'_j,\,E'\rangle$ are lying on it. Their intersection is a line containing the points $B_{ij}$, $B_i$ and $B_j$. So, the point $B_{ij}$ is lying on the line $B_i B_j$ for any $i<j$. Obviously $B_1\notin \mathcal{N}$, and therefore $\dim \mathcal{S}=n-1$ where $\mathcal{S}:=\langle B_1,\,\mathcal{N}\rangle$, and for any $j>1$ we have $B_j\in \mathcal{S}$. So, for any $i$, $j$ such that $1\leq i<j\leq n$ we have $B_{ij}\in \mathcal{S}$, $B_i\in \mathcal{S}$. Thus, $\mathcal{L}_{(\mathcal{R},\,\mathcal{R}')}\subset \mathcal{S}$. Obviously, the opposite inclusion is also hold. $\square$

\smallskip

%%%%%%%%%%%%%%%%%%%%%%%%%%%%%%%%%%%%%%%%%%%%%%%%%%%%%%%%%%%%%%%%%%%%%%%%%%%%%%%%%

\begin{Def}
  We say $\mathcal{L}_{(\mathcal{R},\,\mathcal{R}')}$ to be the Desargues hyperplane generated by $\mathcal{R}$ and $\mathcal{R}'$.
  %({\bf Check the grammar!})
\end{Def}

\smallskip

%=====================================================================
%=====================================================================

\section{Geometrical description of $H$-orbits}

\begin{Th}\label{main result}
  Let $\mathcal{R}$ and $\mathcal{R}'$ be in strict perspective. Then they are equivalent iff the Desargues hyperplane $\mathcal{L}_{(\mathcal{R},\,\mathcal{R}')}$ is passing through $A$.
\end{Th}

\smallskip

{\it Proof}. Let $\mathcal{R}$ and $\mathcal{R}'$ be in strict perspective. Then they are in perspective, and according to Proposition~\ref{Perspectivity} there exist bases $\vec{\mathcal{R}}$ and $\vec{\mathcal{R}}''$ generating these frames such that for some numbers $a_1,\,\ldots,\,a_n,\,h$ ($h\neq0$) the equalities (\ref{eq*4}) hold. Let $(x^0:\,x^1:\,\ldots\,:x^n)$ be the homogeneous coordinates on $\mathbb{P}_n$ with respect to $\mathcal{R}$. In these coordinates hyperplanes $L$ and $L'$ are given by the following equations:
$$
  L\colon\, x^0=0,\quad L'\colon\,a_i x^i-hx^0=0.
$$
Thus, the equation of the bunch $S$ of hyperplanes $S_{(\lambda\colon\mu)}$ passing through $L\cap L'$ one can present as follows
\begin{equation}\label{bunch_S}
  S_{(\lambda\colon\mu)}\colon\,\lambda x^0+\mu a_i x^i=0,\quad \lambda^2+\mu^2\neq0.
\end{equation}
For the point $B_1$ there exists a vector $\vec{B}_1\in\pi^{-1}(B_1)$ such that
$$
  \vec{B}_1=-\frac{e}{a_1}\vec{A}_1+\vec{E},\quad e:=a_1+\ldots+a_n+1-h.
$$
Hyperplane $\mathcal{L}_{(\mathcal{R},\mathcal{R}')}$ is distinguished from the bunch by the condition $B_1\in \mathcal{L}_{(\mathcal{R},\mathcal{R}')}$. The latter imposes the relation on $\lambda$ and $\mu$:
$$
  \lambda=(1-h)\mu.
$$
We substitute this into (\ref{bunch_S}) and obtain the equation of $\mathcal{L}_{(\mathcal{R},\mathcal{R}')}$:
\begin{equation}\label{eq6}
  \mathcal{L}_{(\mathcal{R},\mathcal{R}')}\colon\,(1-h) x^0+ a_i x^i=0,\quad \lambda^2+\mu^2\neq0.
\end{equation}
Therefore
$$
  A_0\in \mathcal{L}_{(\mathcal{R},\mathcal{R}')}\stackrel{(\ref{eq6})}{\Leftrightarrow} h=1\stackrel{Th.~\ref{Criterion_1}}{\Leftrightarrow} \mathcal{R}\sim\mathcal{R}'.\,\square
$$

\smallskip

%=============================================================================
%=============================================================================
\section{Conclusion}
The theorem \ref{main result} shows the existence of surprising relation between the Desargues theorem and the geometry of $H$-orbits.

We can distinguish some applications of the results above:
\begin{itemize}
  \item explicit construction of quotient bundles of the adapted frame bundle over a submanifold $S\subset \mathbb{P}_n$;
  \item geometrical description of linear connections on a submanifold $S\subset \mathbb{P}_n$.
\end{itemize}

One of the further directions of research is studying relations between other classical theorems of projective geometry and the representations group theory.

%%%%%%%%%%%%%%%%%%%%%%%%%%%%%%%%%% Список литературы %%%%%%%%%%%%%%%%%%%%%%%%%%%%%%%%%%%%%

\end{document}